\documentclass{svmult}
\usepackage{graphicx}        
\usepackage{multicol}        
\usepackage[bottom]{footmisc}
\usepackage{graphics}
\usepackage{epsfig}
\usepackage{amssymb}
\usepackage{amsfonts}
\usepackage{amsmath}
\usepackage{mathrsfs}

\begin{document}

\title*{An iterative method for solving non-linear hydromagnetic  equations}
\author{C. BOULBE\inst{1}\and
T.Z. BOULMEZAOUD\inst{2} \and T. AMARI\inst{3}}
\institute{Laboratoire de math\'ematiques appliqu\'ees, \\
Universit\'e de Pau
  et des Pays de l'Adour, \\ IPRA - Av de l'Universit\'e, 64000 PAU,
  FRANCE. \\
\texttt{c.boulbe@etud.univ-pau.fr} \and Laboratoire de
math\'ematiques appliqu\'ees, \\ Universit\'e de Versailles Saint
Quentin, \\ 45 av des Etats Unis,
  75035 Versailles, FRANCE. \\
  \texttt{boulmeza@math.uvsq.fr}\and Centre de Physique
  Th\'eorique,\\ Ecole Polytechnique, F91128 Palaiseau Cedex,FRANCE.
   \\ \texttt{amari@cpht.polytechnique.fr}}

%
%
\maketitle
\begin{abstract}
We propose an iterative finite element method  for  solving
non-linear hydromagnetic and steady Euler's equations. Some
three-dimensional computational tests are given to confirm the
convergence and the  high efficiency of the method.
\end{abstract}

\section{Introduction. Statement of the problem}
\label{sec:1}
The understanding  of plasma equilibria is one of the most important
problems in magnetohydrodynamics and arises in several fields including solar physics and  thermonuclear fusion. Such an equilibria is often governed
by the well known steady  hydromagnetic equations
\begin{eqnarray}
{\mathrm{\vec{curl}}\,}\vec{B}\times\vec{B}+\nabla p=0\label{eqn:MS},\\
{\mathrm{div}\,}\vec{B}=0\label{eqn:div1},
\end{eqnarray}
which describe the balance of the Lorentz force by
 pressure. Here $\vec{B}$ and   $p$ are respectively the magnetic
field and the pressure.

Notice that system
 (\ref{eqn:MS})+(\ref{eqn:div1}) is quite similar to steady
 inviscid fluid equations
\begin{eqnarray}
\vec{v}.\nabla \vec{v}+\nabla p=0,\\ {\mathrm{div}\,}\vec{v}=0.
\end{eqnarray}
This analogy is due to the vectorial identity $$ \displaystyle{\vec{v} . \nabla
\vec{v} - \nabla \frac{|\vec{v}|^2}{2} = {\mathrm{\vec{curl}}\,} \vec{v} \times \vec{v}}. $$  System of
equations (\ref{eqn:MS})+(\ref{eqn:div1}) must be completed
 with some boundary conditions on $\vec{B}$ and $p$. Physical
considerations suggest to prescribe the boundary normal field
component:
\begin{equation}
\vec{B}.\vec{n}=g \mbox{ on }\partial\Omega
\end{equation}
where  $g$ satisfies the compatibility condition
$\displaystyle{\int_{\Omega}g=0 }$ due to the  equation ${\mathrm{div}\,} \vec{B} = 0$.
Defining the   inflow boundary as
$\Gamma^-=\{\vec{x}\in\Omega,\;\;\vec{B}(\vec{x}).\vec{n}(\vec{x})<0\},$ one can also
prescribe  the normal component ${\mathrm{\vec{curl}}\,} \vec{B} . \vec{n}$ of  the current
density and the pressure $p$ on $\Gamma^-$
\begin{equation}
{\mathrm{\vec{curl}}\,}\vec{B}.\vec{n}=h\mbox{ on }\Gamma^-,
\end{equation}
\begin{equation}
p=p_{0} \mbox{ on }\Gamma^-.
\end{equation}
One can notice that if the pressure is neglected, equations
(\ref{eqn:MS})+(\ref{eqn:div1})  become
\begin{eqnarray}
{\mathrm{\vec{curl}}\,}\vec{B}\times\vec{B}=\vec{0}\label{eqn:fff},\\ {\mathrm{div}\,}\vec{B}=0\label{eqn:div2}.
\end{eqnarray}
Equation (\ref{eqn:fff}) means that the magnetic field
and its curl, which represents the current density, are everywhere
aligned. The magnetic field is said  {\it
Beltrami} or {\it force-free} (FF). A usual way to tackle the problem
(\ref{eqn:fff}) + (\ref{eqn:div2}) consists to rewrite
  equation
(\ref{eqn:fff}) into the form
\begin{equation}\label{eqn:fff2}
{\mathrm{\vec{curl}}\,} \vec{B}=\lambda(\vec{x})\vec{B},
\end{equation}
where $\lambda(\vec{x})$ is a scalar function which  can be a constant
function or can depend on $\vec{x}$. In the former, the $\vec{B}$ field is
said linear FF. In the latter, it is said non linear. \\ Some
partial  results concerning existence of 3D solutions of equations
(\ref{eqn:MS})+(\ref{eqn:div1}) in bounded domains are given in
\cite{alber} and \cite{laurence1}. Linear force-free-fields were
studied in \cite{boulmeza_m2an}. For the existence of non-linear
ones the reader can refer to \cite{boulmeza_zamp},
\cite{boulmeza_note2}.

The numerical solving of equations (\ref{eqn:MS})+(\ref{eqn:div1})
and equations (\ref{eqn:fff})+(\ref{eqn:div2}) is of importance in
magnetohydrodynamics studies and in solar physics. As it is known,
the reconstruction of the  coronal magnetic field
  has  is of a  great utility in observational and theoretical
 studies of the magnetic structures in the solar atmosphere. In
 this paper, we propose an iterative  process for solving these
 equations (section \ref{sec:2}).  A finite element method
 is proposed for solving each one of the arising problems.

\section{An iterative method for the magnetostatic system}
\label{sec:2}

Our objective here is to expose an iterative method for solving
the non-linear equations (\ref{eqn:MS})+(\ref{eqn:div1}) in a
bounded and simply-connected domain. The starting idea of the
method consists to split the current density $\vec{\omega}={\mathrm{\vec{curl}}\,} \vec{B}$ into
the sum
\begin{equation}\label{omega123}
\vec{\omega}=\vec{ \omega_{||}} + \vec{ \omega_{\perp}},
\end{equation}
 where the vector field $\vec{ \omega_{||}} = \mu(\vec{x}) \vec{B}$ is
collinear to $\vec{B}$, while $\vec{ \omega_{\perp}}$ is perpendicular to $\vec{B}$. The
problem is decomposed formally  into a
curl-div system on $\vec{B}(\vec{x})$ and two first order hyperbolic
equations on $\mu(\vec{x})$ and $p(\vec{x})$. \\ More precisely, writing $\vec{ \omega_{||}}(\vec{x})=\mu(\vec{x})\vec{B}(\vec{x})$
where $\mu$ is a scalar function and taking the divergence of
(\ref{omega123}), gives
\begin{equation}\label{eqn:mu}
 \vec{B}.\nabla\mu=-{\mathrm{div}\,}\vec{ \omega_{\perp}}.
\end{equation}
Notice that the pressure satisfies a similar equation since
\begin{equation}
\vec{B}.\nabla p=0.
\end{equation}
Equation (\ref{eqn:MS}) becomes
\begin{equation}
\vec{ \omega_{\perp}} \times \vec{B} = - \nabla p,
\end{equation}
which means that $\vec{ \omega_{\perp}}(\vec{x}) = \displaystyle{\frac{1}{|\vec{B}(\vec{x})|^2}\nabla p(\vec{x})\times \vec{B}(\vec{x})}$
if $|\vec{B}(\vec{x})| \ne \vec{0}$.

In consideration of  these remarks, we are going  now to  propose
an iterative process to solve  non-linear systems (\ref{eqn:MS})+(\ref{eqn:div1}). In this  process  the transport equation (\ref{eqn:mu}) is perturbed
by adding an artificial reaction term. Namely,   we construct a
sequence $(\vec{B}^{(n)}, p^{(n)})_{n \ge 0} $ as follows:
\begin{itemize}
\item The starting guess $\vec{B}_0 \in H^1(\Omega)$ is chosen as the irrotational field
associated to $g$  defined by
\begin{equation}\label{eqn:eq0}
{\mathrm{\vec{curl}}\,} \vec{B}_0 = \vec{0} \mbox{ in } \Omega, \; {\mathrm{div}\,} \vec{B}_0 = 0 \mbox{ in }
\Omega \mbox{ and } \vec{B}_0  . \vec{n} = g\mbox{ on } \partial \Omega.
\end{equation}
This is a usual problem which can be reduced to a scalar Neumann
problem since the domain is simply-connected.
\item   For all $n
\ge 0$, $p^{(n)}$ is  solution of  the system
\begin{equation}\label{eqn:eq1bis}
\left\{
\begin{array}{rcl}
 \vec{B}^{(n)}.\nabla p^{(n)} + \eta p^{(n)} &=& \eta p^{(n-1)}  \mbox{ in }\Omega,\\
 p^{(n)} &=& p_0 \mbox{ on } \partial \Omega,
\end{array}
\right.
\end{equation}
where $\eta$ is a small parameter and $p^{(-1)}=0$.
  \item  For all $n
\ge 0$, $\displaystyle{\vec{ \omega_{\perp}}^{(n)}=\frac{1}{|\vec{B}^{(n)}|^2} \nabla p^{(n)} 
\times \vec{B}^{(n)}}$ and  $\vec{ \omega_{||}}^{(n)} = \mu^{(n)} \vec{B}^{(n)}$,
where $\mu^{(n)}$ satisfies
\begin{equation}\label{eqn:pbmu}
  \left\{
  \begin{array}{rcl}
 \vec{B}^{(n)}.\nabla\mu^{(n)} + \epsilon \mu^{(n)}&=&-{\mathrm{div}\,}\vec{ \omega_{\perp}}^{(n)}+\varepsilon\mu^{(n-1)}\mbox{ in
      }\Omega,\\ \mu^{(n)}(\vec{B}^{(n)}.\vec{n})&=&h-\vec{ \omega_{\perp}}^{(n)}.\vec{n} \mbox{ on }\Gamma^-.
 \end{array}
  \right.
  \end{equation}
  Here $\mu^{(-1)}=0$ and $\epsilon$ is a small parameter.
\item  For all $n
\ge 0$, $\vec{B}^{(n+1)} = \vec{B}_0 + \vec{b}^{(n+1)}$, with $\vec{b}^{(n+1)}$ solution of
      $$
      \left\{
      \begin{array}{lrc}
    {\mathrm{\vec{curl}}\,} \vec{b}^{(n+1)}=\vec{\omega}^{(n)} + \nabla q^{(n)} \mbox{ in }\Omega,\\
    {\mathrm{div}\,} \vec{b}^{(n+1)}=0\mbox{ in }\Omega,\\
    \vec{b}^{(n+1)}.\vec{n}=0\mbox{ on }\partial\Omega,
      \end{array}
      \right.
      $$
where $\vec{\omega}^{(n)} = \vec{ \omega_{||}}^{(n)} + \vec{ \omega_{\perp}}^{(n)}$ while $q^{(n)}$ is solution of the Laplace problem
\begin{equation}\label{eqn:lap_cor}
-\Delta q^{(n)} =  {\mathrm{div}\,} \vec{\omega}^{(n)} \mbox{ in } \Omega, \; \mbox{
and } q^{(n)} = 0 \mbox{ on } \partial  \Omega.
\end{equation}
 \end{itemize}
  Notice  that the appearance of the correction term $\nabla q^{(n)}$ is due to the fact that
${\mathrm{div}\,}(\vec{\omega}^{(n)})$ is not zero in general. \\
 The convergence of this iterative process is not an easy matter.  We conjecture that it converges if $h$ is sufficiently small and $|\vec{B}_0 (\vec{x})|\geq c > 0$ in $\Omega$ for some constant $c > 0$. Nevertheless, in the case of linear force-free fields (in that case the algorithm is simplified
 since at each iteration $p^{(n)} = 0$, $\vec{ \omega_{\perp}}^{(n)} = \vec{0}$ and $\mu^{(n)}$ is a fixed real) Boulmezaoud and Amari
\cite{boulmeza_aplin} proved that this process is
super-convergent. The proof of convergence in the general case is not
given and remains an open question.\\ Notice that the same algorithm can be used
for computing linear or non-linear force-free fields which are
solutions of (\ref{eqn:div2})+(\ref{eqn:fff2}), provided that the
computation of the pressure $p^{(n)}$ and the vector field
$\vec{ \omega_{\perp}}^{(n)}$ are dropped.

\section{Finite element discretization}\label{sec:3}
Here  we give a short description of the finite elements methods
we use for solving problems arising in the iterative process
exposed above. Observe first that at each iteration of the
algorithm one should solve two problems:\\ (a) A
reaction-convection  problem of the form: {\it find $u$ solution
of}
\begin{equation}
\label{eqn:hyperbolic}
\left\{
\begin{array}{rcl}
{\mathrm{div}\,}(u\vec{B}) +\sigma u  &=&f \mbox{ in }\Omega,\\ u&=&
h \mbox{ on }  \Gamma^-.
\end{array}
\right.
\end{equation}
(b) A vector potential problem: {\it find the pair $(\vec{b}, q)$
satisfying}
\begin{equation}\label{eqn:systemB}
\left\{
\begin{array}{rcll}
{\mathrm{\vec{curl}}\,} \vec{b} - \nabla q &=& \vec{j} &\mbox{ in }\Omega,\\ {\mathrm{div}\,}\vec{b} &=& 0& \mbox{
in }\Omega,\\ \vec{b}.\vec{n} &=& 0 &\mbox{ on }\partial\Omega,\\ q&=&0 &\mbox{ on
}\partial \Omega.
\end{array}
\right.
\end{equation}
We begin with the approximation of (\ref{eqn:hyperbolic}).

It is well known that the
direct application of a Galerkin finite
elements method to the singularly perturbed   problem   (\ref{eqn:hyperbolic})
may lead to the appearance of  spurious oscillations and instabilities.
In the two last decades, several methods were proposed to remove this drawback (especially in the two  dimensional case). Among these methods, one
can recall the {\it streamline diffusion} method  (see Brookes and
Hughes \cite{brookes}, see also,  e. g.,  Johnson et
{\it al. } \cite{johnson}), the {\it discontinuous} Galerkin
method (see Lesaint \cite{lesaint}) and  {\it bubble functions}  methods (see, e. g., Brezzi et {\it al. } \cite{brezzi2}).  Here we shall use the method of streamline diffusion.

Thus, let us  consider a family of regular triangulations $({\mathscr T}
_h)$ of $\Omega$. The discrete problem we consider is $$ ({\mathscr
P}_h)\left\{
\begin {array}{l}
\mbox{\it Find $u_h \in W_h$ such that} \\ a_h(u_h,
w_h) = \ell_h(w_h), \; \forall w_h \in W_h,
\end{array}
\right. $$ where $$
\begin{array}{rcl}
a_h(u_h, w_h) &=&\displaystyle{ \int_{\Omega}(\vec{B} . \nabla u_h + \sigma u_h).
(w_h+ \delta_h  \vec{B} . \nabla w_h) dx  -\int_{\Gamma^-} u_h  w_h (\vec{B}
. \vec{n}) dx, } \\ \ell_h(w_h) &=&\displaystyle{ \int_{\Omega} f(\vec{x}) (\vec{B} . \nabla
w_h + \delta_h  w_h)  -\int_{\Gamma^-}   \alpha_0   w_h (\vec{B} . \vec{n})
dx. }
\end{array}
$$
 Here $W_h$ stands for the finite elements space
$$ W_h=\{v_h\in H^1(\Omega); \;v_{|K}\in {\mathbb P}_k(K),\;\forall
K\in{\mathscr T}_h\}, $$ where for each $K \in {\mathscr T}_h$,  ${\mathbb P}_k(K)$ denotes
the space
of polynomials of degree less or equal $k$. \\ One can
prove  that  the problem $({{\mathscr P}_h})$ has a
unique solution $u_h\in W_h$ when  $\delta_h  \sigma  < 1$.
Moreover,  if $\delta_h = c h$ for some constant $c$ and if
 $\vec{B} \in L^{\infty}(\Omega)^3 \cap{H({\mathrm{div}\,}; \, \Omega)}$ and $u \in H^{\ell +
1}(\Omega)$ for some $\ell \ge 1$, then
\begin{equation}
(1- \delta_h \sigma) \||u - u_h\|| \le C h^{\ell +1/2}
\|u\|_{H^{\ell + 1}(\Omega)},
\end{equation}
where  $ \||w\||_\Omega^2 =  \delta_h \|\vec{B} . \nabla w\|_{L^2(\Omega)}^2 +
\sigma \| w\|_{L^2(\Omega)}^2 + \||\vec{B}. \vec{n}|^{1/2} w\|^2_{L^2( \partial\Omega)}. $\\
 Now, we deal with the approximation of the curl-div
system (\ref{eqn:systemB}), which can be dispatched into two
problems: a variational problem (${\mathscr Q}$) in terms of $\vec{b}$
and the fictitious unknown $\theta = 0$, and Laplace equation
(\ref{eqn:lap_cor}) in terms of $q$. We only deal with the
approximation of $\vec{b}$, since we shall see that the computation of
the $q$ is useless. 
Denote by ${H({\mathrm{\vec{curl}}\,}; \, \Omega)}$ the space $$ {H({\mathrm{\vec{curl}}\,}; \, \Omega)}=\{\vec{v}\in
L^2(\Omega)^3;\;{\mathrm{\vec{curl}}\,}\vec{v}\in L^2(\Omega)^3\} $$ equipped with its
usual norm. The statement of problem (\ref{eqn:systemB}) suggests
the use of
 an ${H({\mathrm{\vec{curl}}\,}; \, \Omega)}$ approximation. Define $M$ the space
$$
M=\{ v\in H^1(\Omega),\;\;\int_\Omega v\ dx =0\}.
$$
 Let $X_h \subset {H({\mathrm{\vec{curl}}\,}; \, \Omega)}$, $M_h \subset M$
 two finite-dimensional subspaces and set  $$
V_h=\{\vec{v_h}\in X_h;\;(\vec{v_h},\nabla\mu_h)=0,\;\forall\mu_h\in M_h\}. $$
We make the following assumptions
$$
\begin{array}{l}
 ({\mathscr H}_1) \mbox{ the inclusion } \{\nabla\mu_h,\mu_h\in
M_h\}\subset X_h \mbox{ holds, }  \\
  ({\mathscr H}_2) \mbox{ there exists  a constant $C$ such that }
\end{array}
$$
$$\|\vec{v_h}\|_{0,\Omega}\leq
C\|{\mathrm{\vec{curl}}\,}\vec{v_h}\|_{0,\Omega},\; \forall\vec{v_h}\in V_h.
$$
 The discrete version of problem (\ref{eqn:systemB}) writes $$
({\mathscr Q_h}) \left\{
\begin{array}{lrc}
\mbox{{\it Find ($\vec{b_h}$,$\theta_h$) }}\in X_h\times M_h \mbox{ {\it
such as }}\\ \displaystyle{ \forall\vec{v_h}\in X_h, \int_{\Omega}{\mathrm{\vec{curl}}\,}\vec{b_h}.{\mathrm{\vec{curl}}\,}\vec{v_h}
dx +\int_{\Omega}\vec{v_h}.\nabla\theta_h dx=\int_{\Omega}\vec{j}.{\mathrm{\vec{curl}}\,}\vec{v_h} dx,}\\
\displaystyle{\forall\mu_h\in M_h, \int_{\Omega}\vec{b_h}.\nabla\mu_h=0.}
\end{array}
\right. $$ According to Amrouche and {\it al.} \cite{amrouche},
the problem $({\mathscr Q_h})$ has one and only one solution
$(\vec{b}_h,\theta_h)$ with  $\theta_h=0$, and
\begin{equation}\label{cea_hrot}
\|\vec{b}-\vec{b_h}\|_{{H({\mathrm{\vec{curl}}\,}; \, \Omega)}}\le C\inf_{\vec{v_h}\in X_h}\|\vec{b}-\vec{v_h}\|_{{H({\mathrm{\vec{curl}}\,}; \, \Omega)}}.
\end{equation}
A simple manner for constructing the spaces $X_h$ and $M_h$  is to
use  the  $H({\mathrm{\vec{curl}}\,})$ conforming elements of N\'edelec \cite{nedelec}
(see Amrouche and {\it al.}). In that case, the following estimate
 holds
\begin{equation}
\|\vec{b} - \vec{b_h}\|_{H({\mathrm{\vec{curl}}\,}; \Omega)}  \le C h^{\ell} \{ |\vec{b}|_{\ell, \Omega} +
|\vec{b}|_{\ell+1, \Omega}\},
\end{equation}
which is valid if $\vec{b} \in H^{\ell + 1}(\Omega)$. \\
 An important feature of the discrete system $({\mathscr Q_h})$ 
is that only the discrete vector field $\vec{b_h}$ is really unknown. Actually, we know
that   $\theta_h = 0$ . This property can be exploited from a practical viewpoint to
 reduce the discrete system to a smaller one by eliminating $\theta_h$.   In term of matrices, the system writes
\begin{equation}\label{eqn:syst}
\left(
\begin{array}{ccc}
A^{curl}&B^T\\ B&0
\end{array}
\right) \left(
\begin{array}{ccc}
X\\ Y
\end{array}
\right)
=
\left(
\begin{array}{ccc}
C^{curl}\\  0
\end{array}
\right)
\end{equation}
where $A^{curl}$ is a symmetric and positive square matrix ($A^{curl}$ is not
definite neither invertible). We can state the following
\begin{lemma}\label{matrix_lemma}
Let $\Lambda$ be a square positive, definite and symmetric matrix
having the same size as $A$.  Then, the pair $(X, Y)$ is solution
of (\ref{eqn:syst}) if and only if $Y = 0$ and $X$ is solution of
\begin{equation}\label{eqn:dec_syst}
(A^{curl} +  B^T \Lambda B) X = C^{curl}.
\end{equation}
\end{lemma}

\begin{remark}
In Lemma \ref{matrix_lemma}, the matrix $A^{curl}$ and the RHS $C^{curl}$ 
are not arbitrary. Indeed,  if  $G$ denotes  the matrix of
the operator $\nabla:\, M_h\rightarrow X_h$, then necessarily
 $G^T A^{curl} = 0 $ and $G^T C^{curl} = 0$. These identities are the discrete counterpart of the continuous relations ${\mathrm{div}\,} ({\mathrm{\vec{curl}}\,} .) = 0$ and ${\mathrm{div}\,}\vec{j}=0$.
\end{remark}
 A serious  advantage of the new system
(\ref{eqn:dec_syst}) comparing with (\ref{eqn:syst}) is that number of unknowns is reduced.

\section{Computational tests}
\label{sec:4} 
In this last section, we expose some computational
results we obtain with a 3D code. This code use the iterative
method and the finite elements discretization exposed above to
solve problem (\ref{eqn:MS})+(\ref{eqn:div1}) and problem
(\ref{eqn:fff2})+(\ref{eqn:div2}). We compare the exact solution
and the numerical solution and we show the behavior of the errors
in terms of $h$. Two exact solutions are used for the tests.
\begin{itemize}

\item {\bf Test 1 (a non-linear force-free-field). }\\
Let  $(r,\theta,z)$ the cylindrical coordinates with respect to a
point $(x_0,y_0,0)$ ($x_0=-3$ and $y_0=-3$). The pair $(\vec{B}, p)$ is
given $\displaystyle{\vec{B}=\frac{1}{\sqrt{r}}(\vec{e}_{\theta}+
\vec{e}_z)}, \; p(\vec{x}) =
0. $ This is a non-linear force-free field with
$\displaystyle{\lambda=\frac{1}{2r}}.$  Table \ref{tab:1} shows the
behavior of the residue $\|\vec{B}^{(n+1)} - \vec{B}^{(n)}\|_{0, \Omega}$ and
the product ${\mathrm{\vec{curl}}\,} \vec{B}^{(n)}\times\vec{B}^{(n)}$ versus the iteration
number. This example illustrates the superconvergence of the
algorithm.

\begin{table}
\centering \caption{Evolution of $\frac{\|\vec{B}^{(n+1)} - \vec{B}^{(n)}\|_{0,
\Omega}}{\|\vec{B}^{(n)}\|_{0, \Omega}}$ and $\|{\mathrm{\vec{curl}}\,}
\vec{B}^{(n)}\times\vec{B}^{(n)}\|_{\infty}$.} \label{tab:1}
\begin{tabular}{lll}
\hline\noalign{\smallskip} $n$  &  $\frac{\|\vec{B}^{(n+1)} -
\vec{B}^{(n)}\|_{0, \Omega}}{\|\vec{B}^{(n)}\|_{0, \Omega}}$ &  $\|{\mathrm{\vec{curl}}\,}
\vec{B}^{(n)}\times\vec{B}^{(n)}\|_{\infty}$
\\ \noalign{\smallskip}\hline\noalign{\smallskip} 0 & 0.09912&   6.740e-15
\\ 1 & 0.00566& 0.06781  \\ 2 & 0.00036& 0.01939  \\ 3 & 2.644e-05&  0.01910 \\
\noalign{\smallskip}\hline
\end{tabular}
\end{table}

\item {\bf Test 2: (Bennet pinch) }. The pair $(\vec{B}, p)$ is given
by $$\displaystyle{\vec{B}=\nabla A\times \vec{e}_y}\mbox{ and }
\displaystyle{p=\frac{\lambda}{2}\E^{2A}} \mbox{ with } \displaystyle{A=-\ \ln(\frac{1+\lambda k^2 (x^2+z^2)}{2k})}.$$ In table \ref{tab:2},
the relative $L^2$ errors on $\vec{B}$ and $p$ after convergence of the
algorithm are shown. These error decreases as
 $h^{1.8}$, which confirms the high  accuracy of the method.

\begin{table}
\centering \caption{Relative errors on $\vec{B}_h$ and $p_h$ in norm
$L^2$ (test 2).}
\label{tab:2}       
%
%
\begin{tabular}{lll}
\hline\noalign{\smallskip} $h$ &
$\frac{\|\vec{B}-\vec{B}_h\|_{0,\Omega}}{\|\vec{B}\|_{0,\Omega}}$ &
$\frac{\|p-p_h\|_{0,\Omega}}{\|p\|_{0,\Omega}}$  \\
\noalign{\smallskip}\hline\noalign{\smallskip} 0.69282 & 0.03837
& 0.08648 \\ 0.23094 & 0.00492 & 0.01102 \\ 0.13856 &
0.00191 & 0.00396 \\ 0.09897 & 0.00108 & 0.00201\\
\noalign{\smallskip}\hline
\end{tabular}
\end{table}

\begin{figure}
\centering
   \begin{minipage}[c]{.46\linewidth}
      \includegraphics[width=6.1cm]{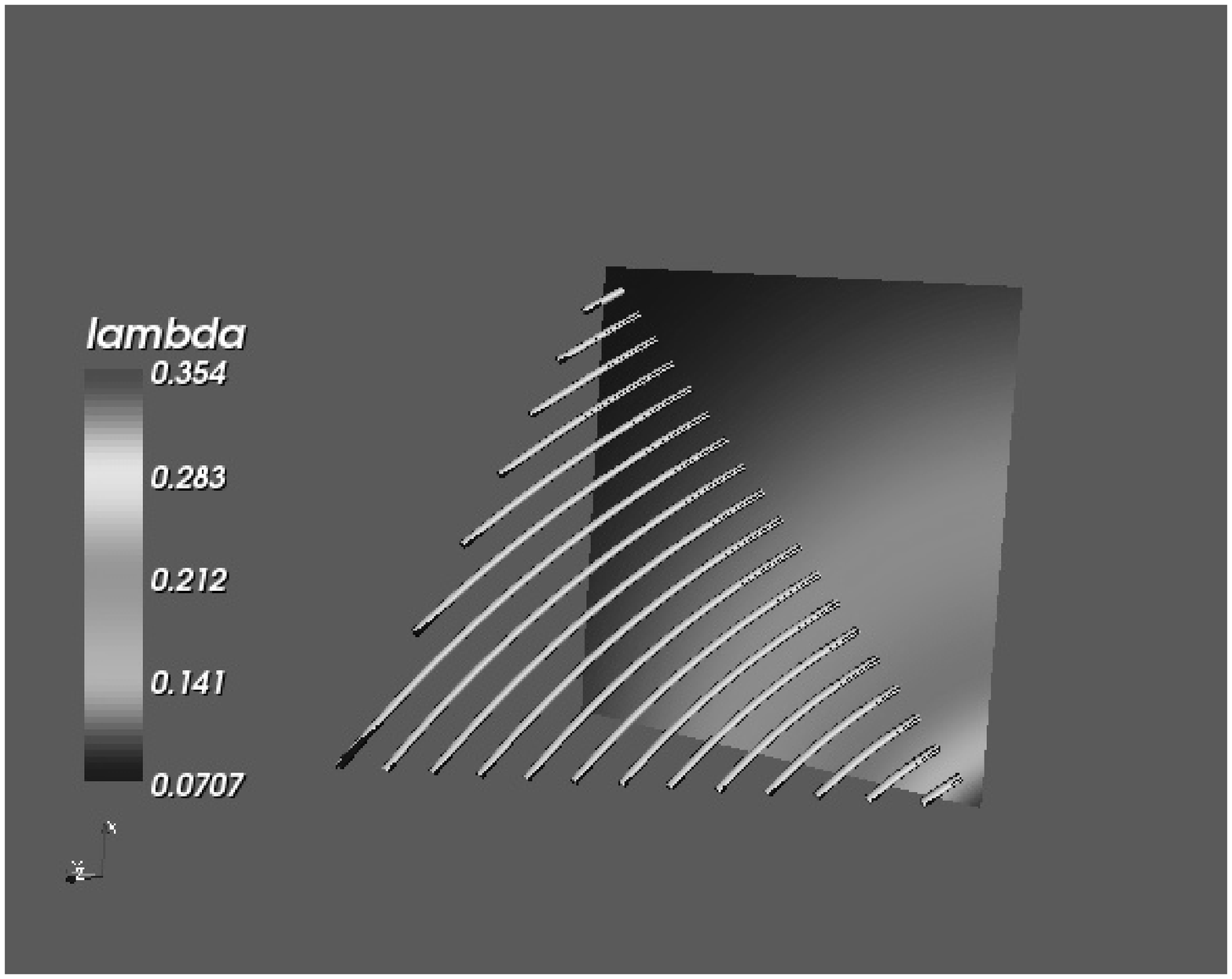}
   \end{minipage} \hfill
   \begin{minipage}[c]{.46\linewidth}
      \includegraphics[width=6.1cm]{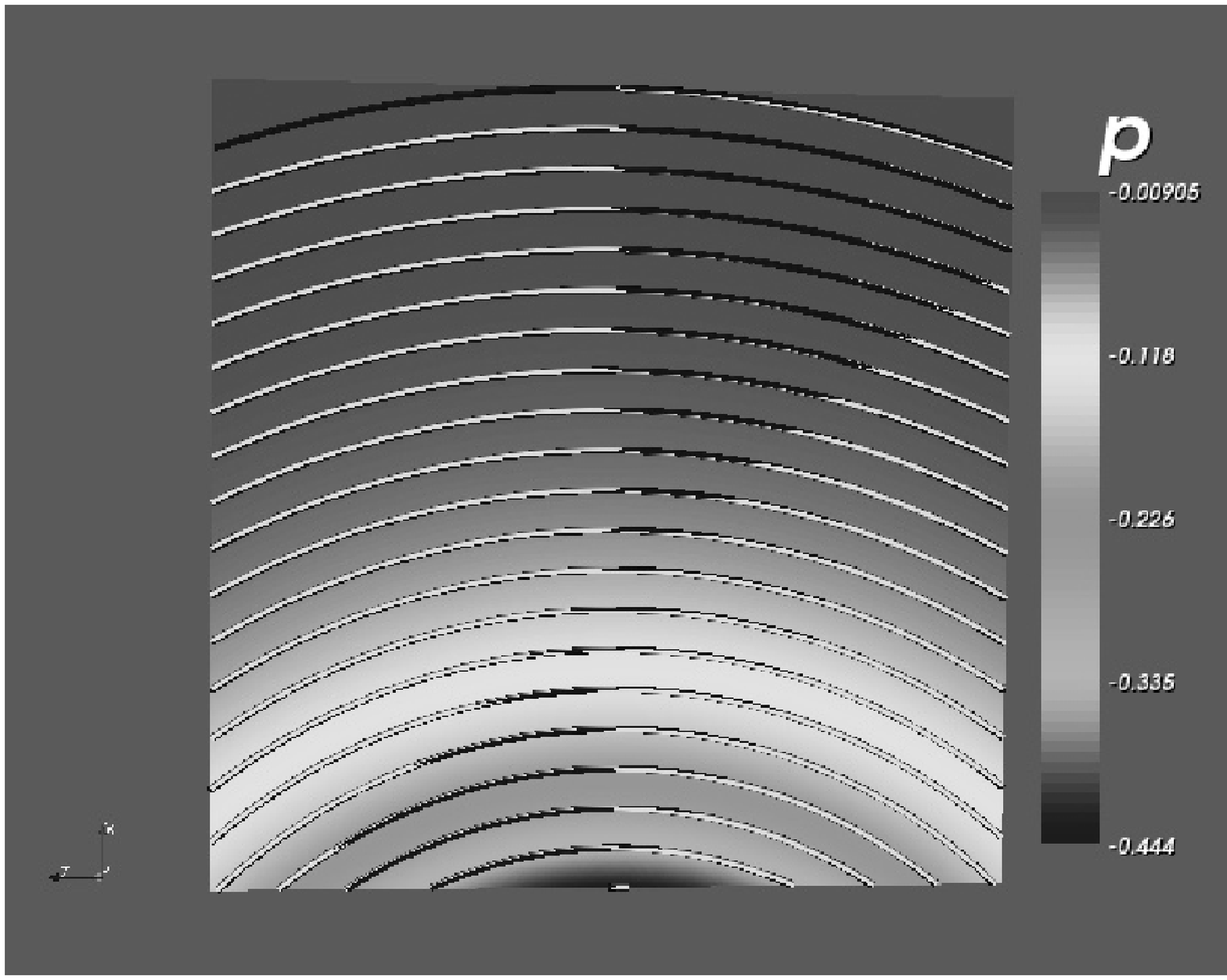}
   \end{minipage}
\caption{Superposition of the the exact and the numerical
  solutions in the case of test 1 on the left and  in a $(x-z)$ plane
  2D cut for the test 2 on the right.}
\end{figure}

\end{itemize}
%
%




\begin{thebibliography}{}
%
%
\bibitem{alber}  Alber, H.D.: Existence of three-dimensional,
 steady, inviscid, incompressible flows with non-vanishing vorticity.
  Math. Ann.  \textbf{292}, 493--528 (1992)
\bibitem{amrouche} Amrouche, C., Bernardi, C ., {D}auge, M.,
   Girault, V.: Vector potentials in three-dimensional
 non-smooth
    domains, Math. Methods Appl. Sci. \textbf{9}, 823--864 (1998).


\bibitem{boulmeza_note2} Boulmezaoud, T.Z.:  On the existence of
  non-linear Beltrami fields.  Comptes Rendus de l'Acad\'emie des
  Sciences, \textbf{328},  437--442 (1999).

\bibitem{boulmeza_m2an} Boulmezaoud, T.Z.,  Maday, Y., Amari, T.:
 On the linear Beltrami fields in bounded and unbounded
three-dimensional domains: Mathematical Modelling and Numerical
Analysis, \textbf{33}, 359--394 (1999)

\bibitem{boulmeza_zamp} Boulmezaoud, T.Z. , Amari, T.:
On the existence of non-linear force-free fields
       in three-dimensional multiply-connected domains. \`a
 para\^{\i}tre dans  Zeitschrift f\"{u}r Angewandte Mathematik und Physik (ZAMP).

\bibitem{boulmeza_aplin} Boulmezaoud, T. Z., Amari, T.: Approximation
  of linear force-free fields in bounded 3-D
  domains. Math. Comput. Modelling  \textbf{31}, 109--129  (2000)

\bibitem{brezzi2}Brezzi, F.,  Hauke, G., Marini, L.D.,  Sangalli, G.:
  Link-Cutting Bubbles for the Stabilization of Convection-Diffusion-Reaction
Problems, Math. Models Methods Appl. Sci., \textbf{13}, 445--461 (2003)

\bibitem{brookes}  Brooks, A., Hughes,T. R. J.:  Streamline
 upwind/Petrov-Galerkin formulations for convection dominated flows
 with particular emphasis on the incompressible Navier-Stokes
 equations
. FENOMECH '81, Part I (Stuttgart, 1981). Comput. Methods
 Appl. Mech. Engrg.  \textbf{32},  199--259 (1982).

\bibitem{girault} Girault, V., Raviart, P. A.:  Finite element methods
  for Navier-Stokes equations. Springer-Verlag (1986).

\bibitem{johnson} Johnson, C., N\"avert, U., Pitk\"aranta, J.:  Finite element methods for linear hyperbolic
problems, Comp. Meth. Appl. Mech. Engin. \textbf{45},  285--312 (1984).

\bibitem{laurence1} Laurence, P., Bruno, 0.: Existence of 3D toroidal
  MHD equilibria with nonconstant pressure, Communications on Pure and
  Applied Math,  \textbf{49}, 717--764  (1996)

\bibitem{lesaint}Lesaint, P.:  Sur la r\'esolution des syst\`emes  hyperboliques du premier ordre par des m\'ethodes d'\'el\'ements finis
, Th\`ese de Doctorat, UPMC, Paris (1975).

\bibitem{nedelec} N\'ed\'elec, J.C.:  Mixed finite elements in $ R\sp{3}$.
  Numer. Math.  \textbf{35}, 315--341 (1980).



\end{thebibliography}
\end{document}